\documentclass[11pt, reqno]{amsart}

\author[P.~Leonetti]{Paolo Leonetti}
\address{Department of Economics, Universit\'a degli Studi dell'Insubria, via Monte Generoso 71, 21100 Varese, Italy}
\email{leonetti.paolo@gmail.com}
\urladdr{\url{https://sites.google.com/site/leonettipaolo/}} 


\title{Rough families, cluster points, and cores}

\usepackage[T1]{fontenc}
\usepackage{amsmath}
\usepackage{amssymb}
\usepackage{amsthm}
\usepackage[left=3.5cm, right=3.5cm, paperheight=11.8in]{geometry}
\usepackage{hyperref}
\usepackage{fancyhdr}
\usepackage{enumitem}
\usepackage{comment}
\usepackage{nicefrac}
\usepackage{mathrsfs}
\usepackage{bm}
\usepackage{graphicx}
\usepackage[utf8]{inputenc}
\usepackage{cancel}
\usepackage{mathtools}

\AtBeginDocument{%
   \def\MR#1{}
}

\newtheorem{thm}{Theorem}[section]
\newtheorem{cor}[thm]{Corollary}

\newtheorem{prop}[thm]{Proposition}

\theoremstyle{definition} 
\newtheorem{defi}[thm]{Definition}
\let\olddefi\defi
\renewcommand{\defi}{\olddefi\normalfont}
\newtheorem{example}[thm]{Example}
\let\oldexample\example
\renewcommand{\example}{\oldexample\normalfont}
\newtheorem{rmk}[thm]{Remark}
\let\oldrmk\rmk
\renewcommand{\rmk}{\oldrmk\normalfont}

\hypersetup{
    pdftitle={A note on rough ideal convergence},
    pdfauthor={Paolo Leonetti},
    pdfmenubar=false,
    pdffitwindow=true,
    pdfstartview=FitH,
    colorlinks=true,
    linkcolor=blue,
    citecolor=green,
    urlcolor=cyan
}

\providecommand{\MR}[1]{}

\providecommand{\MR}{\relax\ifhmode\unskip\space\fi MR }

\providecommand{\href}[2]{#2}

\begin{document}

\maketitle
\thispagestyle{empty}

\begin{abstract}
We define the notion of ideal convergence for sequences $(x_n)$ with values in topological spaces $X$ with respect to a 
family $\{F_\eta: \eta \in X\}$ of subsets of $X$ with $\eta \in F_\eta$. 
Each set $F_\eta$ 
quantifies the degree of accuracy of the convergence toward $\eta$. 

After proving that this is really a new notion, we provide some properties of the set of limit points and characterize the latter through the ideal cluster points and the ideal core of $(x_n)$. 
\end{abstract}


\section{Introduction and Main Results}\label{sec:intro}

Let $\mathcal{I}\subseteq \mathcal{P}(\omega)$ be an ideal on the nonnegative integers $\omega$, that is, a family closed under subsets and finite unions. 
It is also assumed that the family of finite subsets of $\omega$, denoted by $\mathrm{Fin}$, is contained in $\mathcal{I}$ and that 
and $\omega \notin \mathcal{I}$. 
%

Let also $\bm{x}=(x_n)$ be a sequence taking values in a topological space $(X,\tau)$ (note that it is not assumed to be Hausdorff). Lastly, let 
$$
\mathscr{F}:=\{F_\eta: \eta \in X\}
$$
be a \textbf{rough family}, that is, a collection of subsets of $X$ with the property that $\eta \in F_\eta$ for all $\eta \in X$. Rough families, as it will be clear from the following definition, quantifies the \textquotedblleft degree of accuracy\textquotedblright\, of sequences taking values in $X$ toward their limits $\eta$. In particular, they can change depending on $\eta$: smaller sets $F_\eta$ can be interpreted as smaller oscillations of the tail of sequence around its limit $\eta$. 
\begin{defi}\label{defi:mainroughness}
A sequence $\bm{x}=(x_n)$ is said to be $\mathcal{I}$\emph{-convergent to} $\eta \in X$ \emph{with roughness} 
$\mathscr{F}$, shortened as 
$
(\mathcal{I}, \mathscr{F}, \tau)\text{-}\lim\nolimits_n x_n =\eta,
$ 
provided that 
$$
\{n \in \omega: x_n \notin U\} \in \mathcal{I}
$$ 
for all $\tau$-open sets $U\subseteq X$ such that $F_\eta \subseteq U$. We denote by $\mathsf{L}_{\bm{x}}(\mathcal{I}, \mathscr{F}, \tau)$ the set of all 
$\mathcal{I}$-limits of $\bm{x}$ with roughness $\mathscr{F}$, that is, $$
\mathsf{L}_{\bm{x}}(\mathcal{I}, \mathscr{F}, \tau):=\left\{\eta \in X: (\mathcal{I}, \mathscr{F}, \tau)\text{-}\lim\nolimits_n x_n =\eta\right\}.
$$
\end{defi}

Notice that: 
\begin{enumerate}[label={\rm (\roman{*})}]
\item if $F_\eta=X$ for all $\eta \in X$, then $\mathcal{I}$-convergence with roughness $\mathscr{F}$ corresponds to ordinary convergence with respect to the trivial topology $\tau_0:=\{\emptyset, X\}$;
\item if $F_\eta=\{\eta\}$ for all $\eta \in X$, then $\mathcal{I}$-convergence with roughness $\mathscr{F}$ simplifies to the classical $\mathcal{I}$-convergence with respect to the same topology $\tau$; in such case, we simply speak about $(\mathcal{I}, \tau)$-convergence, 
see e.g. \cite{MR3920799}; 
\item if $F_\eta=\{\eta\}$ for all $\eta \in X$ and, in addition, $\mathcal{I}=\mathrm{Fin}$, then $\mathcal{I}$-convergence with roughness $\mathscr{F}$ corresponds to ordinary $\tau$-convergence; 
\item special instances where $X$ is a normed vector space and each $F_\eta$ is chosen as the closed ball with center $\eta$ and fixed radius $r \in [0,\infty)$ have been studied in several works, see e.g. \cite{MR2412408, MR2412409, MR1990651} and references therein. 
\end{enumerate}



It is remarkable that Definition \ref{defi:mainroughness} may \emph{not} correspond to $(\mathcal{J},\nu)$-convergence, for every ideal $\mathcal{J}$ on $\omega$ and for every topology $\nu$ on $X$:
\begin{prop}\label{prop:necessaryprop}
Suppose 
that $X=\mathbb{R}$ is endowed with the standard Euclidean topology $\tau$. Then there exists a rough family $\mathscr{F}$ such that, for each ideal $\mathcal{I}$ on $\omega$, there is no ideal $\mathcal{J}$ on $\omega$ and no topology $\nu$ on $\mathbb{R}$ for which the equivalence
\begin{equation}\label{eq:equivalencecontradiction}
(\mathcal{I}, \mathscr{F}, \tau)\text{-}\lim\nolimits_n x_n=\eta 
\quad \text{ if and only if } \quad 
(\mathcal{J}, \nu)\text{-}\lim\nolimits_n x_n=\eta 
\end{equation}
holds for all real sequences $(x_n)$ and all $\eta \in \mathbb{R}$. 
\end{prop}
This proves that the type of convergence stated in Definition \ref{defi:mainroughness} defines a new notion which is not included in the classical one. Note that such preliminary result is necessary to avoid unnecessary repetitions of known facts, as it already happened in the literature with other variants of ideal convergence, see for instance the case of \textquotedblleft ideal statistical convergence\textquotedblright\, in \cite[Theorem 2.3]{MR4052262}. 
Hereafter, the dependence on the underlying topology $\tau$ will be made implicit whenever it is clear from the context, so that we will simply write $(\mathcal{I}, \mathscr{F})\text{-}\lim_n x_n=\eta$ or $\mathsf{L}_{\bm{x}}(\mathcal{I}, \mathscr{F})$. 

The aim of this note is to prove some characterizations of $\mathcal{I}$-convergence with roughness $\mathscr{F}$. 
For, we need to recall some definitions. A point $\eta \in X$ is said to be an $\bm{\mathcal{I}}$\textbf{-cluster point} of a sequence $\bm{x}$ if $\{n \in \omega: x_n \in U\} \notin \mathcal{I}$ for all open sets $U$ containing $\eta$. The set of $\mathcal{I}$-cluster points of $\bm{x}$ is denoted by $\Gamma_{\bm{x}}(\mathcal{I})$. It is known that $\Gamma_{\bm{x}}(\mathcal{I})$ is a closed subset of $X$, and it is nonempty provided that $\{n \in \omega: x_n \notin K\} \in \mathcal{I}$ for some compact $K\subseteq X$. Moreover, it follows readily from the definitions that 
$$
\mathrm{L}_{\bm{x}}(\mathcal{I}, \mathscr{F}) \subseteq \Gamma_{\bm{x}}(\mathcal{I}). 
$$
We refer to \cite{MR3920799} for basic properties and characterizations of $\mathcal{I}$-cluster points. 
%


\begin{thm}\label{thm:charact}
Let $\bm{x}$ be a sequence taking values in a regular topological space $X$ such that $\{n \in \omega: x_n \notin K\} \in \mathcal{I}$ for some compact set $K\subseteq X$. 
Also, let $\mathcal{I}$ be an ideal on $\omega$, let $\mathscr{F}$ be a rough family, and pick $\eta \in X$ such that $F_\eta$ is closed. 
Then 
$$
(\mathcal{I}, \mathscr{F})\text{-}\lim\nolimits_n x_n=\eta 
\quad \,\,\text{ if and only if }\quad \,\,
\Gamma_{\bm{x}}(\mathcal{I})\subseteq F_\eta.
$$
\end{thm}

Note that the hypothesis on $\bm{x}$ includes the case of relatively compact sequences (which corresponds to the case $\mathcal{I}=\mathrm{Fin}$). In addition, the claim does not hold without any restriction of $F_\eta$: for, suppose that $X=\mathbb{R}$, $F_\eta=(\eta-\nicefrac{1}{2}, \eta+\frac{1}{2})$ for all $\eta \in \mathbb{R}$, $\mathcal{I}=\mathrm{Fin}$ and $\bm{x}$ is an enumeration of the rationals in $[0,1]$. Then it is readily checked that $\mathrm{L}_{\bm{x}}(\mathcal{I}, \mathscr{F})=\{\nicefrac{1}{2}\}$ and, on the other hand, there are no $\eta$ for which $[0,1]=\Gamma_{\bm{x}}(\mathcal{I})\subseteq F_\eta$. 

The following corollary is immediate: 
\begin{cor}\label{cor:firstrepresent}
Suppose, in addition to the hypotheses of Theorem \ref{thm:charact}, that every $F_\eta$ is closed. Then 
$$
\mathsf{L}_{\bm{x}}(\mathcal{I}, \mathscr{F})=\left\{\eta \in X: \Gamma_{\bm{x}}(\mathcal{I})\subseteq F_\eta\right\}. 
$$
\end{cor}

%
%

Hereafter, if $X$ is a metric space with metric $d$, we denote the closed ball with center $\eta \in X$ and radius $r \in [0,\infty]$ by 
$$
B_r(\eta):=\{x \in X: d(x,\eta)\le r\}. 
$$
In particular, $B_0(\eta)=\{\eta\}$ and $B_\infty(\eta)=X$.

As a [non-]linear property of $(\mathcal{I}, \mathscr{F})$-convergence, we obtain the following:
\begin{prop}\label{prop:consequencenotvectorspace}
Let $X$ be a normed vector space, let $\mathcal{I}$ be a nonmaximal ideal on $\omega$, and fix a rough family $\mathscr{F}$ for which the sets $F_\eta$ are uniformly bounded. 
Then the family of $(\mathcal{I}, \mathscr{F})$-convergent sequences is a vector space if and only if 
$F_\eta=\{\eta\}$ for all $\eta \in X$.
\end{prop}

We remark that, if $\mathcal{I}$ is maximal (that is, if its dual filter $\mathcal{I}^\star:=\{S\subseteq \omega: \omega\setminus S \in \mathcal{I}\}$ is a free ultrafilter on $\omega$), then all relatively compact sequences are $\mathcal{I}$-convergent (hence also $(\mathcal{I}, \mathscr{F})$-convergent, for each rough family $\mathscr{F}$).



Given a topological space $X$, we endow the hyperspace 
$$
\mathcal{H}(X):=\{F\subseteq X: F \text{ nonempty closed}\}.
$$
with the \emph{upper Vietoris topology} $\widehat{\tau}$, that is, the topology generated by the base of sets $\{F \in \mathcal{H}(X): F\subseteq U\}$, with $U \in \tau$ open. 
Moreover, we recall that a metric space $X$ is said to have the $\textsc{UC}$\emph{-property} if nonempty closed sets are at a positive distance apart, 
that is, for all $F, F^{\prime} \in \mathcal{H}(X)$ with 
$F\cap F^\prime=\emptyset$, there exists $\varepsilon>0$ such that $d(x,x^\prime)> \varepsilon$ for all $x \in F$ and $x^\prime \in F^\prime$, where $d$ is the metric on $X$. See  \cite{MR1165059, MR42109} and references therein. (It is remarkable that a metric space $X$ has the $\textsc{UC}$-property if and only if the ordinary Vietoris topology is weaker than the Hausdorff topology on $\mathcal{H}(X)$. In addition, standard Euclidean spaces $\mathbb{R}^k$ have the $\textsc{UC}$-property.)

\begin{thm}\label{thm:limitsetclosed}
Let $\bm{x}$ be a sequence taking values in a topological space $X$, let $\mathcal{I}$ be an ideal on $\omega$, and pick a rough family $\mathscr{F}$ made by closed sets. Also, suppose that the map 
$\eta\mapsto F_\eta$ is $\widehat{\tau}$-continuous. 
Then $\mathrm{L}_{\bm{x}}(\mathcal{I}, \mathscr{F})$ is closed. 
\end{thm}

The result above does not hold, similarly, without any restriction on the rough family $\mathscr{F}$. 
Indeed, suppose that $X=\mathbb{R}$, $F_\eta=(\eta-3,\eta+3)$ for all $\eta \in \mathbb{R}$, $\mathcal{I}=\mathrm{Fin}$, and that $\bm{x}$ is defined by $x_n=(-1)^n$ for all $n \in \omega$. Then $\mathrm{L}_{\bm{x}}(\mathcal{I}, \mathscr{F})=(-2,2)$. In particular, together with the example given after Theorem \ref{thm:charact}, if every $F_\eta$ is open then $\mathrm{L}_{\bm{x}}(\mathcal{I}, \mathscr{F})$ is not necessarily closed, nor open. 

In the special case where $X$ is a metric space with the $\textsc{UC}$-property and each $F_\eta$ is a closed ball 
$B_{r(\eta)}(\eta)$, for some function $r(\cdot)$,  
we obtain the following: 
\begin{cor}\label{cor:UCproperty}
Let $\bm{x}$ be a sequence taking values in a metric space $X$ with the \textsc{UC}-property, let $\mathcal{I}$ be an ideal on $\omega$, and fix an upper semicontinuous function $r: X\to [0,\infty)$ such that $F_\eta=B_{r(\eta)}(\eta)$ for all $\eta \in X$.  
Then $\mathrm{L}_{\bm{x}}(\mathcal{I}, \mathscr{F})$ is closed. 
\end{cor}

%

In the case that $X$ has a linear structure on $X$, we can show that $\mathrm{L}_{\bm{x}}(\mathcal{I}, \mathscr{F})$ is convex: 
\begin{thm}\label{thm:convex}
Let $\bm{x}$ be a sequence taking values in a normed vector space $X$ with the \textsc{UC}-property, let $\mathcal{I}$ be an ideal on $\omega$, and fix a concave function $r: X\to [0,\infty)$ such that $F_\eta=B_{r(\eta)}(\eta)$ for all $\eta \in X$.  
Then $\mathrm{L}_{\bm{x}}(\mathcal{I}, \mathscr{F})$ is convex. 
\end{thm}

%


Using the above results, we provide a relationship between $(\mathcal{I}, \mathscr{F})$-convergence and the $\bm{\mathcal{I}}$\textbf{-core} of a sequence $\bm{x}$, see \cite{MR4126774, MR3955010}. For, given an ideal $\mathcal{I}$ on $\omega$ and a sequence $\bm{x}$ taking values in a topological vector space $X$, we define 
$$
\mathrm{core}_{\bm{x}}(\mathcal{I}):=\bigcap_{E \in \mathcal{I}^\star} \overline{\mathrm{co}}(\{x_n: n \in E\}).
$$
In other words, the $\mathcal{I}$-core of $\bm{x}$ is the least closed convex set containing $\{x_n: n \in E\}$ for all $E \in \mathcal{I}^\star$ (where 
$\overline{\mathrm{co}}$ stands for the closed convex hull operator). In the case where $\mathcal{I}=\mathrm{Fin}$, we obtain the so-called \emph{Knopp core}, see \cite{MR2241135, MR1441452, MR0142952} and references therein.

\begin{thm}\label{thm:coreconseq}
Let $\bm{x}$ be a sequence taking values in a locally convex space $X$ such that $\{n \in \omega: x_n \notin K\} \in \mathcal{I}$ for some compact $K\subseteq X$. Also, let $\mathcal{I}$ be an ideal on $\omega$ and pick a rough family $\mathscr{F}$ such that every $F_\eta$ is closed and convex. Then
\begin{equation}\label{eq:claimedequality}
\mathrm{L}_{\bm{x}}(\mathcal{I}, \mathscr{F})=\left\{\eta \in X: \mathrm{core}_{\bm{x}}(\mathcal{I})\subseteq F_\eta\right\}.
\end{equation}
\end{thm}

We remark that the hypothesis on $\bm{x}$ cannot be removed. Indeed, if $X=\ell_\infty$ is the Banach space of bounded real sequences, endowed with the supremum norm, $F_\eta=B_1(\eta)$ for all $\eta \in \ell_\infty$ and $\bm{x}=(e_0,-e_0,e_1,-e_1,\ldots)$, where $e_k$ stands for the $k$th unit vector (so that $\bm{x}$ is not relatively compact), then it is readily seen that
$$
\mathrm{core}_{\bm{x}}(\mathrm{Fin})=\{0\}
\quad \text{ and }\quad 
\mathrm{L}_{\bm{x}}(\mathrm{Fin}, \mathscr{F})=c_{00}.
$$
(Here, $c_{00}$ represents the Banach subspace of eventually zero sequences.) However, $(1,1,\ldots) \in \{x \in \ell_\infty: 0 \in B_1(x)\}\setminus c_{00}$. 
 To sum up, $\bm{x}$ is a nonconvergent bounded sequence, its Knopp core 
is a singleton, it is $(\mathrm{Fin}, \mathscr{F})$-convergent to every sequence $\eta \in c_{00}$, and the claimed equality \eqref{eq:claimedequality} fails. 

Lastly, we prove that the above properties coincide for relatively compact sequences:
\begin{cor}\label{cor:equivalenceIcore}
Let $\bm{x}$ be a sequence taking values in a metric vector space $X$ such that $\{n \in \omega: x_n \notin K\} \in \mathcal{I}$ for some compact $K\subseteq X$. Also, let $\mathcal{I}$ be an ideal on $\omega$, fix $r \in [0,\infty)$, and define $\mathscr{F}$ by $F_\eta=B_r(\eta)$ for all $\eta \in X$. 

Then the following are equivalent\textup{:}
\begin{enumerate}[label={\rm (\roman{*})}]
\item \label{eq:condition1} $\bm{x}$ is $\mathcal{I}$-convergent to some $\eta \in X$\textup{;}
\item \label{eq:condition2} $\mathrm{core}_{\bm{x}}(\mathcal{I})$ contains a unique vector $\eta^\prime \in X$\textup{;}
\item \label{eq:condition3} $\mathrm{L}_{\bm{x}}(\mathrm{Fin}, \mathscr{F})=B_r(\eta^{\prime\prime})$ for some $\eta^{\prime\prime} \in X$.
\end{enumerate}
In addition, in such case, $\eta=\eta^\prime=\eta^{\prime\prime}$. 
\end{cor}

The proofs follow in the next section. 

\section{Proofs} 

\begin{proof}
[Proof of Proposition \ref{prop:necessaryprop}]
Let $\mathscr{F}$ be the rough family defined by $F_\eta:=[\eta-1,\eta+1]$ for all $\eta \in \mathbb{R}$. Fix also an ideal $\mathcal{J}$ on $\omega$ and a topology $\nu$ on $\mathbb{R}$, and suppose for the sake of contradiction that equivalence \eqref{eq:equivalencecontradiction} holds for all $\bm{x}$ and $\eta$. We divide the proof in two cases. 

\medskip

\textsc{Case (i): $\mathcal{I}$ not maximal.} If $\mathcal{I}$ is not maximal, there exists a partition $\{A,B\}$ of $\omega$ such that $A,B \notin \mathcal{I}$. 
For each $r \in \mathbb{R}$ and $h \in (0,1]$, let $\bm{x}^{(r,h)}$ be the real sequence defined by $x_n^{(r)}=r$ if $n \in A$ and $x_n^{(r)}=r+h$ otherwise. 
Define similarly $\bm{y}^{(r,h)}$ such that $y_n^{(r)}=r+h$ if $n \in A$ and $x_n^{(r,h)}=r$ otherwise. 
If follows that 
$$
\mathsf{L}_{\bm{x}^{(r,h)}}(\mathcal{I}, \mathscr{F}, \tau)=
\mathsf{L}_{\bm{y}^{(r,h)}}(\mathcal{I}, \mathscr{F}, \tau)=
[r+h-1,r+1]
$$ 
for each $r \in \mathbb{R}$ and $h \in (0,1]$. Now, let $V$ be a nonempty $\nu$-open set and fix a point $r \in V$. Fix also $h \in (0,1]$. By the equivalence \eqref{eq:equivalencecontradiction}, we obtain that $
\{n \in \omega: x^{(r,h)}_n\notin V\} \in \mathcal{J}
$ and $
\{n \in\omega: y^{(r,h)}_n \notin V\} \in \mathcal{J}.
$
If $r+h\notin V$, this can be rewritten as $B \in \mathcal{J}$ and $A \in \mathcal{J}$, respectively, which is impossible because it would imply $\omega=A\cup B \in \mathcal{J}$. Hence $r+h \in V$. By the arbitrariness of $h$, we obtain $[r,r+1]\subseteq V$. However, since $r$ is arbitrary, we conclude that $V=\mathbb{R}$, therefore $\nu$ stands for the trivial topology $\tau_0$. This is a contradiction because $\mathsf{L}_{\bm{x}^{(r,h)}}(\mathcal{I}, \mathscr{F}, \tau)\neq \mathbb{R}$. 

\medskip 

\textsc{Case (ii): $\mathcal{I}$ maximal.} If $\mathcal{I}$ is maximal, then either $A:=\{2n: n\in \omega\}$ or $B:=\{2n+1: n \in \omega\}$ belong to $\mathcal{I}$. Suppose without loss of generality $B \in \mathcal{I}$ (the remaining case is symmetric). With the same notations above, it follows that 
$$
\mathsf{L}_{\bm{x}^{(r,h)}}(\mathcal{I}, \mathscr{F}, \tau)=[r-1,r+1] 
\quad \text{ and }\quad 
\mathsf{L}_{\bm{y}^{(r,h)}}(\mathcal{I}, \mathscr{F}, \tau)=
[r+h-1,r+h+1]
$$
for each $r \in \mathbb{R}$ and $h \in (0,1]$. Since $r \in \mathsf{L}_{\bm{x}^{(r,h)}}(\mathcal{I}, \mathscr{F}, \tau) \cap \mathsf{L}_{\bm{y}^{(r,h)}}(\mathcal{I}, \mathscr{F}, \tau)$, we obtain by the equivalence \eqref{eq:equivalencecontradiction} that $(\mathcal{J}, \nu)\text{-}\lim_n x_n^{(r,h)}=(\mathcal{J}, \nu)\text{-}\lim_n y_n^{(r,h)}=r$. Hence the sequence $\bm{x}^{(r,h)}+\bm{y}^{(r,h)}$, which is constantly equal to $2r+h$, is $(\mathcal{J},\nu)$-convergent to $2r$. Since $r$ and $h$ are arbitrary, it follows that $\nu$ is the trivial topology $\tau_0$, reaching the same contradiction as in the nonmaximal case above. 
\end{proof}

\begin{proof}
[Proof of Theorem \ref{thm:charact}]
\textsc{If part.} 
Suppose that every $\mathcal{I}$-cluster point belongs to $F_\eta$. Fix an open set $U\subseteq X$ containing $\Gamma_{\bm{x}}(\mathcal{I})$. We need to show that $S:=\{n \in \omega: x_n \notin U\}$ belongs to $\mathcal{I}$. 
For, notice that $\Gamma_{\bm{x}}(\mathcal{I})$ is a nonempty compact subset of $K$, see \cite[Lemma 3.1]{MR3920799}. Suppose by contradiction that $S\notin \mathcal{I}$. Considering that $\{n \in \omega: x_n\in K\setminus U\}=S\setminus I \notin \mathcal{I}$ and that $K\setminus U$ is compact, we conclude, again by \cite[Lemma 3.1]{MR3920799}, that $\Gamma_{\bm{x}}(\mathcal{I}) \cap (K\setminus U)\neq \emptyset$, which contradicts the hypothesis $S\notin \mathcal{I}$. Therefore $(\mathcal{I}, \mathscr{F})\text{-}\lim_n x_n=\eta$. 


\medskip

\textsc{Only If part.} 
Suppose for the sake of contradiction that there exists an $\mathcal{I}$-cluster point $\eta_0 \in X\setminus F_\eta$. Since $X$ is regular, one can pick disjoint open sets $U_0, U_\eta\subseteq X$ such $\eta_0 \in U_0$ and $F_\eta \subseteq U_\eta$.  However, this  is impossible because 
$\{n \in \omega: x_n \in U_0\}$, which does not belong to $\mathcal{I}$ since $\eta_0$ is an $\mathcal{I}$-cluster point of $\bm{x}$, is contained in $\{n \in \omega: x_n \notin U_\eta\}$, which belongs to $\mathcal{I}$ because $(\mathcal{I}, \mathscr{F})\text{-}\lim\nolimits_n x_n=\eta$. 
\end{proof}

\begin{rmk}\label{rmk:ifandonlyifassumptions}
It is clear from the proof above that the \textsc{If part} of Theorem \ref{thm:charact} holds for arbitrary topological spaces $X$ and arbitrary rough families $\mathscr{F}$. 

On the same line, the \textsc{Only If part} holds also for arbitrary sequences. 
\end{rmk}

\begin{proof}
[Proof of Proposition \ref{prop:consequencenotvectorspace}] 
%
\textsc{If part} This is a folklore fact, by the linearity of the $\mathcal{I}$-limit. 

\medskip

\textsc{Only If part.} Since $\mathcal{I}$ is not maximal, there exists a partition $\{A,B\}$ of $\omega$ such that $A,B \notin \mathcal{I}$. Suppose also that there exists $\eta \in X$ such that $F_\eta \neq \{\eta\}$, hence it is possible to fix a point $\eta^\prime \in F_\eta\setminus \{\eta\}$. Let $\bm{x}$ be the sequence such that $x_n=\eta$ if $\eta \in A$ and $x_n=\eta^\prime$ otherwise. It follows by the definition of $(\mathcal{I}, \mathscr{F})$-convergence that $\{\eta, \eta^\prime\}\subseteq \mathrm{L}_{\bm{x}}(\mathcal{I}, \mathscr{F})$. Pick also $r \in (0,\infty)$ such that $\mathrm{diam}(F_\eta) \le r$ for all $\eta \in X$. If the claim were false, the sequence $k\bm{x}$ would be $(\mathcal{I}, \mathscr{F})$-convergent for all $k \in \omega$. Notice that $\Gamma_{k\bm{x}}(\mathcal{I})=k\Gamma_{\bm{x}}(\mathcal{I})=\{k\eta, k\eta^\prime\}$, see e.g.  \cite[Proposition 3.2]{MR4428911}. 
However, the distance between the latter two $\mathcal{I}$-cluster points can be made arbitrarily large as $k\to \infty$, which contradicts the hypothesis that the sets $F_\eta$ are uniformly bounded. 
\end{proof}

\begin{proof}[Proof of Theorem \ref{thm:limitsetclosed}]
If $\mathrm{L}_{\bm{x}}(\mathcal{I}, \mathscr{F})=\emptyset$, the claim is obvious. Otherwise, pick a $\tau$-convergent net $(\eta_i)_{i \in I}$ with values in $\mathrm{L}_{\bm{x}}(\mathcal{I}, \mathscr{F})$ and define $\eta:=\lim_i \eta_i$. Since the map $\eta\mapsto F_\eta$ is $\widehat{\tau}$-continuous, then the net $(F_{\eta_i})_{i \in I}$ is $\widehat{\tau}$-convergent to $F_\eta$. Fix an arbitrary open set $U\subseteq X$ which contains $F_\eta$ and define the $\widehat{\tau}$-open set $\widehat{U}:=\{F \in \mathcal{H}(X): F\subseteq U\}$. By the convergence of $(F_{\eta_i})_{i \in I}$, there exists $j \in I$ such that $F_{\eta_j} \in \widehat{U}$, i.e., $F_{\eta_j}\subseteq U$. Since $\eta_j \in \mathrm{L}_{\bm{x}}(\mathcal{I}, \mathscr{F})$, it follows that $\{n \in \omega: x_n\notin U\} \in \mathcal{I}$. Therefore $(\mathcal{I}, \mathscr{F})\text{-}\lim_nx_n=\eta$.
\end{proof}

\begin{proof}
[Proof of Corollary \ref{cor:UCproperty}] 
Denote by $d$ the metric on $X$. Thanks to Theorem \ref{thm:limitsetclosed}, it is sufficient to show that the map $\eta\mapsto F_\eta$ is $\widehat{\tau}$-continuous. Pick a convergent net $(\eta_i)_{i \in I}$ with limit $\eta \in X$, hence $\lim_i d(\eta_i, \eta)=0$. We claim that the net of closed balls $(B_{r(\eta_i)}(\eta_i))_{i \in I}$ is $\widehat{\tau}$-convergent to $B_{r(\eta)}(\eta)$. For, pick an open set $U$ containing $B_{r(\eta)}(\eta)$. In particular, $\eta \in U$. 
As in the previous proof, set $\widehat{U}:=\{F \in \mathcal{H}(X): F\subseteq U\}$. If $U=X$ then $\{n \in \omega: x_n \notin U\}=\emptyset \in \mathcal{I}$. Otherwise $U$ is a proper subset of $X$, hence $X\setminus U$ is a nonempty closed set disjoint from $F_\eta$. 
Since $X$ has the $\mathrm{UC}$-property, it follows that there exists $\varepsilon>0$ such that $d(x,y)\ge \varepsilon$ for all $x \in F_\eta$ and $y \in X\setminus U$. At this point, set
$$
G:=\left\{x \in X: d(x,\eta)<r(\eta)+\varepsilon\,\right\}. 
$$
It follows by contruction that $F_\eta \subseteq G\subseteq U$. By the upper semicontinuity of $r$ and the convergence of $(\eta_i)_{i \in I}$, there exists an index $i_0 \in I$ such that $r(\eta_i)<r(\eta)+\nicefrac{\varepsilon}{2}$ and $d(\eta_i,\eta)<\nicefrac{\varepsilon}{2}$ for all $i\ge i_0$. Therefore
$$
\forall i\ge i_0, \forall x \in F_{\eta_{i_0}}, \quad d(x,\eta)\le d(x,\eta_i)+d(\eta_i,\eta) \le r(\eta_j)+\nicefrac{\varepsilon}{2} < r(\eta)+\varepsilon. 
$$
This shows that $F_{\eta_i}\subseteq G\subseteq U$ (hence, $F_{\eta_j} \in \widehat{U}$) for all $i \ge i_0$, concluding the proof. 
\end{proof}

\begin{proof}
[Proof of Theorem \ref{thm:convex}]
If $\mathrm{L}_{\bm{x}}(\mathcal{I}, \mathscr{F})=\emptyset$, the claim is obvious. Otherwise, fix two vectors $\eta, \eta^\prime \in \mathrm{L}_{\bm{x}}(\mathcal{I}, \mathscr{F})$, a weight $\alpha \in (0,1)$, and define $\gamma:=\alpha \eta+(1-\alpha) \eta^\prime$. We claim that $\gamma \in \mathrm{L}_{\bm{x}}(\mathcal{I}, \mathscr{F})$. For, pick an open set $U$ containing $F_\gamma$. If $U=X$ then $\{n \in \omega: x_n \notin U\}=\emptyset \in \mathcal{I}$. If $U\neq X$, then by the \textsc{UC}-property of $X$, there exists $\varepsilon>0$ such that $F_\gamma\subseteq V\subseteq U$, where $V$ is the open ball with center $\gamma$ and radius $r(\gamma)+\varepsilon$. Lastly, set
$$
S:=\left\{n \in \omega: \|x_n-\eta\|\ge r(\eta)+\varepsilon \text{ or } \|x_n-\eta^\prime\|\ge r(\eta^\prime)+\varepsilon\right\}.
$$
Since $\eta, \eta^\prime \in \mathrm{L}_{\bm{x}}(\mathcal{I}, \mathscr{F})$, it follows that $S \in \mathcal{I}$. At this point, for each $n \in \omega\setminus S$, 
\begin{displaymath}
\begin{split}
\|x_n-\gamma\|&=\|\alpha (x_n-\eta)-(1-\alpha)(x_n-\eta^\prime)\|\\
&\le \alpha \|x_n-\eta\|+(1-\alpha)\|x_n-\eta^\prime\|\\
&< \alpha (r(\eta)+\varepsilon)+(1-\alpha)(r(\eta^\prime)+\varepsilon) 
\le r(\gamma)+\varepsilon.
\end{split}
\end{displaymath}
Therefore $\{n \in \omega: x_n \notin U\}\subseteq \{n \in \omega: x_n \notin V\}\subseteq S\in \mathcal{I}$. Since $U$ is arbitrary, we conclude that $\gamma \in \mathrm{L}_{\bm{x}}(\mathcal{I}, \mathscr{F})$. 
\end{proof}

\begin{rmk}\label{rmk:generalizationmetricugly}
It is worth noting that the analogue of Theorem \ref{thm:convex} holds in metrizable vector spaces $X$ with a compatible metric $d$ which is translation invariant and for which $d(\alpha x,0) \le \alpha d(x,0)$ for all $x \in X$ and $\alpha \in (0,1)$. 
\end{rmk}

\begin{proof}
[Proof of Theorem \ref{thm:coreconseq}]
Since every topological vector space is regular, 
it follows by Corollary \ref{cor:firstrepresent} and the hypothesis that each $F_\eta$ is closed and convex that
$$
\mathrm{L}_{\bm{x}}(\mathcal{I}, \mathscr{F})=\left\{\eta \in X: \overline{\mathrm{co}}(\Gamma_{\bm{x}}(\mathcal{I}))\subseteq F_\eta\right\}.
$$
The conclusion follows by \cite[Theorem 2.2]{MR3955010} and \cite[Theorem 3.4]{MR4126774}, which state that $\mathrm{core}_{\bm{x}}(\mathcal{I})$ coincides with $\overline{\mathrm{co}}(\Gamma_{\bm{x}}(\mathcal{I}))$. 
\end{proof}

\begin{proof}
[Proof of Corollary \ref{cor:equivalenceIcore}]
\ref{eq:condition1} $\Longleftrightarrow$ \ref{eq:condition2}. See \cite[Proposition 3.2]{MR3955010}. 

\medskip

\ref{eq:condition1} $\Longrightarrow$ \ref{eq:condition3}. Suppose that $\mathcal{I}\text{-}\lim_n x_n=\eta$. Pick $\gamma \in B_r(\eta)$ and an open set $U$ containing $B_r(\gamma)$. Since $\eta \in U$, it follows by $\{n \in\omega: x_n \notin U\} \in \mathcal{I}$. Therefore $\bm{x}$ is $(\mathcal{I}, \mathscr{F})$-convergent to $\gamma$. Conversely, if $\gamma \notin B_r(\eta)$, then $\eta \notin B_r(\gamma)$. Since $X$ is regular, there exists disjont open sets $U, U_\gamma$ such that $\eta \in U$ and $B_r(\gamma)\subseteq U_\gamma$. However, this implies that $\{n \in \omega: x_n \notin U_\gamma\}\supseteq \{n \in \omega: x_n \in U\} \in \mathcal{I}^\star$. By the arbitrariness of $\gamma$ in both cases, we conclude that $\mathrm{L}_{\bm{x}}(\mathrm{Fin}, \mathscr{F})=B_r(\eta)$. 

\medskip

\ref{eq:condition3} $\Longrightarrow$ \ref{eq:condition2}. Thanks to Theorem \ref{thm:coreconseq}, we obtain necessarily that $\mathrm{core}_{\bm{x}}(\mathcal{I})=\{\eta\}$. 
\end{proof}

\bibliographystyle{amsplain}
\bibliography{ideale}

\end{document}